\documentclass[12pt,a4paper]{amsart}

\usepackage{natbib}
\usepackage[margin=3cm]{geometry}

\usepackage{mathtools}
\usepackage{booktabs}
\usepackage{tikz}
\usetikzlibrary{cd} 

\newcommand{\invset}{\mathcal{I}}
\newcommand{\tij}[1]{t_{#1}}
\newcommand{\n}{\mathbf n}
\newcommand{\id}{e}
\newcommand{\rws}{\mathcal{R}}
\newcommand{\cayley}{\mathcal{C}}
\newcommand{\writesto}{\rightarrow}
\newcommand{\writestoc}[1]{\underset{\mathclap{#1}}{\rightarrow}}

\newcommand{\N}{\mathbb N}
\newcommand{\newl}{\ell}

\numberwithin{equation}{section}
\newtheorem{prop}{Proposition}[section]

\newtheorem{definition}[prop]{Definition}

\newtheorem{theorem}[prop]{Theorem}

\begin{document}

\title[Weighted genome rearrangement distance]{A path-deformation framework for determining weighted genome rearrangement distance}

\author[Bhatia et al]{Sangeeta Bhatia\,$^{1,3}$, Attila Egri-Nagy\,$^{1,4}$, Stuart Serdoz\,$^{1}$, Cheryl E. Praeger\,$^{2}$, Volker Gebhardt\,$^{1}$, Andrew Francis\,$^{1,*}$}
\address{$^{1}$Centre for Research in Mathematics and Data Science, Western Sydney University, Sydney, NSW, Australia \\
$^{2}$School of Physics, Mathematics, and Computing, University of Western Australia, Perth, WA, Australia\\
$^{3}$current address: MRC Centre for Global Infectious Disease Analysis, School of Public Health, Imperial College London, UK\\
$^{4}$current address: Akita International University, Japan\\ 
$^*$Corresponding author.} 

\begin{abstract}
Measuring the distance between two bacterial genomes under the inversion process is usually done by assuming all inversions to occur with equal probability.  Recently, an approach to calculating inversion distance using group theory was introduced, and is effective for the model in which only very short inversions occur.  In this paper, we show how to use the group-theoretic framework to establish minimal distance for any weighting on the set of inversions, generalizing previous approaches.  To do this we use the theory of rewriting systems for groups, and exploit the Knuth--Bendix algorithm, the first time this theory has been introduced into genome rearrangement problems.  

The central idea of the approach is to use existing group theoretic methods to find an initial path between two genomes in genome space (for instance using only short inversions), and then to deform this path to optimality using a confluent system of rewriting rules generated by the Knuth--Bendix algorithm.
\end{abstract}
\date{\today}
\maketitle

\section{Introduction}
\label{sec:intro}

Large scale changes in the arrangement of genes within a chromosome abound in biology and are key agents of sequence evolution \citep{beckmann2007copy,belda2005genome}. The differences in the order of genes along a chromosome were used as a phylogenetic marker as early as 1938~\citep{dobzhansky1938inversions} when Dobzhansky used them to determine different strains of \emph{Drosophila melanogaster}. \emph{Inversions} of chromosomal fragments are believed to be the main type of rearrangement event in bacterial genomes \citep{belda2005genome}.  

The first formalization of the problem of determining rearrangement distance between gene arrangements was done by \citet{watterson1982chromosome}. A number of methods have been proposed since then to determine the distance between genome arrangements in terms of a single rearrangement operator or a combination of rearrangement operators. In addition to inversions, researchers have considered translocations of chromosomal fragments \citep{bafna1995sorting,bafna1998sorting,yin2013sorting}, fission/fusion of chromosomes, duplication of sequences \citep{chaudhuri2006tandem}, deletion/insertion and a combination of these different operators \citep{yancopoulos2005efficient}. These choices of allowable operations constitute \emph{models of rearrangement}, in which the genomes in the data are assumed to only change according to specific rearrangement operators being considered.

The rearrangement distance between a pair of genomes is usually defined as the minimal number of events from the set of allowed operations required to transform one of the genomes in the pair into the other. For instance, in determining inversion distance between two genomes, the set of legal operations consists of all possible inversions on a gene sequence.  Initial solutions in the case of inversions involve finding the smallest number of inversion events between two genomes and the distance was the count of the events. Thus, each inversion event carries the same weight. If the weight assigned to an inversion event represents the probability of that event, then a model where all events have the same weight can be thought of as finding distances under the uniform distribution.  This model is used in the~\citet{hannenhalli1999transforming} approach, which draws a graph based on the genomes and calculates the minimal distance as a function of features of the graph (for example the number of cycles and paths). 
These methods are simple and fast and have been implemented in software for use by the research community \citep{sankoff1991derange,tesler2002grimm,shao2014exact}. 

As pointed out above, an implicit assumption underlying most of these methods is that all rearrangement operators included in the model are equally probable and thus are given the same weight in the rearrangement distance. 
An inversion model that lies at the other extreme is one that allows only very short inversions. 
A group-theoretic model for sorting circular permutations using inversions acting on two adjacent regions was described by ~\citet{egrinagy2013group}. In a similar vein, \citet{galvao2015sorting} presented an approximation algorithm for sorting signed permutations by only length 2 reversals while \citet{chen1996sorting} gave a characterisation of linear and circular permutations that can be sorted by only length $k$ reversals, for a fixed $k$. 

The biological evidence, however, points somewhere between these two extremes. For example, focusing only on the evidence related to inversions, several studies have suggested that inversions of a short chromosomal fragment are more frequent than that of longer fragments \citep{seoighe2000prevalence, lefebvre2003detection,Darling2008,eisen2000evidence}. Similarly, \citet{seoighe2000prevalence} found a high prevalence of short inversions in the yeast genome. They observed that the conservation of a small neighbourhood of genes, without absolute conservation of order or orientation, suggests that small DNA inversions have contributed significantly to the evolution of ascomycete genomes. In an analysis of four pairs of related bacterial genomes, \citet{lefebvre2003detection} report an over-representation of short inversions, especially those involving a single gene, in comparison with a random inversion model. Analysis of the genome of \emph{Y. Pestis} has also found that all inversions were shorter than expected under a neutral model \citep{Darling2008}.

In view of this information, a natural extension to the definition of rearrangement distance 
that allows for assigning weights (derived from empirical information) to the rearrangement operators, and calculates the minimal \emph{weighted} distance between genome arrangements, might be a better approximation of the underlying biology. 
Thus, it would be useful to have a method to determine weighted inversion distance, where the use of an inversion operator can be penalised based on the number of regions it affects, or where the different operators in a model may be weighted based on type. 

In fact, one of the first algorithms for determining rearrangement distance, proposed in \citet{sankoff1992edit} and \citet{sankoff1992gene}, is in principle capable of assigning weights for inversions and transpositions. An approximation algorithm for sorting a permutation under a particular class of length sensitive cost models, where the cost function is additive i.e., $f(x) + f(y) = f(x+y)$ was presented in \citet{pinter2002genomic}. This approach has been generalized  to a wider class of cost functions \citep{bender2008improved}. This work also improved the bounds on the cost for sorting using an additive cost function. Further pursuing this line of inquiry, \citet{swidan2004sorting} extended the results for signed permutations as well as circular permutations.

In this paper, we present a flexible group-theoretic framework that can be used to determine the weighted rearrangement distance for any model of genome rearrangement in which the rearrangements allowed are invertible.  Thus the framework we propose is applicable to models involving inversions and translocations, but not, for instance, insertions and deletions.  The present work is based on the group theoretic approach of~\citet{egrinagy2013group} and~\citet{francis2014algebraic}. Throughout the paper, we will focus on determining the minimal weighted reversal distance.

\subsection*{Overview of the framework introduced in this paper}

The central idea of the method we propose in this paper is `path deformation' in genome space.  The genome space is the collection of all possible genomes. A path in the genome space is a sequence of genomes where consecutive elements are connected through a single rearrangement operator, and the weight of a path is the sum of the weights of the operators along the path. The minimal weighted distance between two genomes is then the minimal weight of all paths between the genomes. 

To find the minimal weighted distance between two genomes, we start by constructing a path between them.  At the same time, we have also constructed a library of rules in this space. These rules consist of alternate paths, or shortcuts, for a number of small paths. We scan the existing path for any subpaths that could be replaced by a shortcut from our library, generating a new, shorter path. In this way, the existing path is deformed into a new path which is of lower weight than the original path (although it might still not be the least weighted path). Successive iterations of the deformation step should ideally lead us to an optimal path (this is guaranteed only in certain circumstances described below).

The library of ``rewriting rules'' in itself is easy to generate, given a group defined by a presentation (generators and relations, defined in Section~\ref{sec:grppresent}).  The relations, together with the weighting functions, can be transformed to give a set of rewriting rules.  It is also not too difficult to construct an initial path between the two genomes which can be edited using the rules in the library, at least for some models of genome rearrangement. However it is not clear at the outset in what sequence to apply the rules in such a way that one is guaranteed to end with a minimal weight path from one genome to another. This is where the theory of rewriting systems is used.  

A rewriting system that is guaranteed to produce a minimal expression, regardless of the sequence in which the rewriting rules are applied, is called a \emph{confluent} rewriting system.  In this paper, we use the \emph{Knuth-Bendix algorithm}  to transform our rewriting system into a confluent system and use it to construct a minimal weighted path between two genomes given an existing path between them~\citep{knuth1983simple}.

The Knuth--Bendix algorithm is a heuristic whose termination may be affected by the ordering of the generators, which is not an intrinsic property of the input (a group presentation), but rather a choice made when applying the method.  Thus, the input determines neither the result nor the running time of the algorithm, which means its complexity is not defined.  On the other hand, a confluent rewriting system, once obtained, provides a simple algorithm that quickly finds a globally minimal weighted distance between two genomes.  

Note that while this process obtains the global minimal weighted distance, and indeed a path that realises this distance, the path itself is not necessarily unique: several distinct paths through the genome space may attain the globally minimal weighted distance (see~\cite{clark2019bacterial} for a detailed discussion of this).  

We begin the paper by introducing the group based inversion model, and formalising the notion of weighted distance in Section~\ref{subs:invsys}. As a preliminary to the discussion about rewriting systems, we briefly discuss group presentations in Section~\ref{sec:grppresent} and Cayley graphs (Section~\ref{sec:cg}). This is followed by a discussion of rewriting systems and their properties in Section~\ref{sec:rws}. In Section~\ref{sec:examples}, with two concrete examples of rewriting systems, we use weighted distances to draw some phlyogenies. We close off with a discussion of the strengths and limitations of the present work and some directions for future research (Section~\ref{sec:future}).

\section[Inversion systems]{Group theoretic inversion systems}\label{subs:invsys}

The notion of an \emph{inversion system} was formalised in \citet{egrinagy2013group}. Since our work uses much of this language, we briefly summarise the key concepts in this section followed by an extension to a weighted inversion system.  

\medskip
\subsection*{Genomes as permutations and inversion as an action}

A chromosome is represented as a map from a set of positions $\n=\{ 1,2,\dots,n\}$ to a set of regions $X$, usually also labeled with the integers $\n=\{1,2,\dots,n\}$. If we denote the chromosome map by $\pi$, we can write the arrangement in two-line notation as:
\[
\pi=\begin{pmatrix}
1 & 2 & \cdots & n\\
\pi_1 & \pi_2 & \cdots & \pi_n
\end{pmatrix}.
\]
where $\pi_i$ is the region in the $i$'th position. The top row in this view represents the $n$ positions on the chromosome and the bottom row represents the set of regions.

An unsigned inversion operator $\tij{i,j}$ (with $1\le i<j\le n$) in this paradigm is a map from positions to positions. 
When the genome is modeled as a map from positions to regions and a rearrangement operator is a bijection on the set of positions, we require that the rearrangement operator act first on a position and then we map the new position to a region using the genome map,
and so the function composition is from left to right. For a detailed discussion of right and left actions see \citet{bhatia2018position}. The inversion operator $\tij{i,j}$ maps $\pi$ as follows:
\[
\tij{i,j}\pi=
\tij{i,j} \begin{pmatrix}
\cdots & i     & i+1       & \cdots & j     & \cdots\\
\cdots & \pi_i & \pi_{i+1} & \cdots & \pi_j & \cdots
\end{pmatrix} =  
\begin{pmatrix}
\cdots & i & i+1 & \cdots & j & \cdots\\
\cdots & \pi_j & \pi_{j-1} & \cdots & \pi_i & \cdots
\end{pmatrix}
\]
Thus the inversion operator $\tij{i,j}$ flipping regions in positions $i$ to $j$ can be written in cycle notation as follows:
\[ \tij{i,j} := 
\begin{cases}
(i,j)(i+1, j-1)\dots (\frac{i+j}{2} - 1, \frac{i+j}{2} + 1) & \text{ if } j -i \text{ is even,} \\
(i,j)(i+1, j-1)\dots (\frac{i+j-1}{2}, \frac{i+j+1}{2}) & \text{ if } j -i \text{ is odd.} \\
\end{cases} \]
For example, $t_{1,4}=(1,4)(2,3)$, $t_{1,5}=(1, 5)(2, 4)$, and $t_{1,6} = (1, 6)(2, 5)(3, 4)$.

Given genomes $\pi$ and $\pi'$, and a sequence of $k$ inversion operations $t_{i_1,j_1},\dots,t_{i_k,j_k}$ that transform $\pi$ into $\pi'$ when applied in order with $t_{i_1,j_1}$ first, we write
\[\tij{i_k, j_k} \cdots \tij{i_1, j_1}\pi = \pi'.\]
Since $\pi$ is a bijective map from the set of positions to regions, $\pi^{-1}$ is well-defined and we can compose with $\pi^{-1}$, to give
\[\tij{i_k, j_k} \cdots \tij{i_1, j_1}= \pi' \pi^{-1}.\]
Now $\pi' \pi^{-1}$ is a bijective function from positions to positions and therefore an element of the symmetric group on $n$ objects, $S_n$. Thus, the problem of determining a sequence of inversion operations that transforms $\pi$ into $\pi'$ is equivalent to the problem of expressing the group element $\pi'\pi^{-1}$ as a product of the group elements corresponding to the rearrangement operators.

\subsection*{Inversion systems}

An \emph{inversion system} is defined as a tuple $(G, \invset)$ where $G$ is the group of permutations and $\invset$ is a set of inversions such that $G$ is ``generated by'' $\invset$, written $G = \langle \invset \rangle$. In other words, every permutation in $G$ is expressible as a product of elements of $\invset$. 

In general, if we have a subset $S \subset G$ of non-trivial elements from G, then a word over $S$ is a finite sequence of elements of $S$. In this paper, we will assume that $S$ is closed under the operation of taking inverses i.e., for all $s \in S$, we have $s^{-1} \in S$.  We use $S^*$ to represent the (infinite) set of all words over $S$. 

If $S$ generates $G$, then there is a natural map $\Gamma: S^* \rightarrow G$ that sends a word $w = [s_1, s_2,  \cdots ,s_k]$ to the group element $g = s_1 s_2 \cdots s_k$.
The brackets in $w$ are used to emphasise that a word is an ordered sequence of elements of $S$ and to distinguish the sequence from the product $s_1s_2 \cdots s_k$. The set $S^*$ also contains the empty sequence which maps to the identity element of $G$. The length of a group element $g$ with respect to the generating set $S$ is the smallest $r \in \mathbb{N}$ such that there is some element $w \in S^*$, say $w = [s_1, \hdots ,s_r]$, such that $\Gamma(w) = g$, that is, 
\[ s_1 s_2 \cdots s_r = g.\]

The inversion distance between permutations $\pi_1$ and $\pi_2$ is the length of the group element $\pi_2 \pi_1^{-1}$ in the inversion system $(G, \invset)$. For details of inversion systems, the reader is referred to \citet{egrinagy2013group}. 

\subsection*{Weighted length}
The notion of the length of a group element can be extended to the \emph{weighted} length of a group element. Suppose the elements of $S$ are assigned (positive) weights. The weighted length of a word $w = [s_1, s_2,  \dots ,s_k]$ in $S^*$ is the sum of the weights of the $s_i$ where $i$ runs from $1$ through $k$. The weighted length of a group element $g$ is obtained by taking an infimum over the set of all words in $S^*$ that map to $g$.

\begin{definition}[Weighted length]
\label{def:wwl}
Let $S$ be a set of generators of a group $G$. Let $\omega$ be a bounded function $\omega:S \rightarrow \mathbb{R}^+$. The weighted length of a (non-identity) group element $g \in G$ is defined as 
\[
\newl_{S,\omega}(g):= 
\inf \left\{ \sum_{i=1}^{t}{\omega(s_i)}\ \Big\vert\  s_1 s_2  \cdots s_t = g, s_i \in S,t\in\N\right\}.
\]
The weighted length of the identity element $e$ of $G$ is $0$. 
\end{definition}

We define a \emph{weighted inversion system} to be a 3-tuple $(G, \invset, \omega)$ where $G = \langle \invset \rangle$ as before and $\omega : \invset \rightarrow \mathbb{R}^+$.  

\section{Group presentations}\label{sec:grppresent}

We will make use of the important notion of group presentations, from group theory. A \emph{group presentation} is an abstract description of a group $G$ in terms of a generating set $S$ and set of relations $\mathcal{R}$ among the generators. Following \citet[Chapter 1]{coxeter1980generators}, these are defined as follows.

\begin{definition}[Group Presentation]
Let $G$ be a group and let $\id$ be the identity element of $G$.
A \emph{presentation} $\langle S \mid \mathcal{R} \rangle$ for $G$ consists of a generating set $S \subseteq G$ and a set of words $\mathcal{R} \subseteq S^*$ such that
\[ \Gamma(R_i) = \id \quad \text{ for all } R_i \in \mathcal{R}, \] 
and for $w \in S^*$, if $\Gamma(w)  = \id$ then $w$ is an \emph{algebraic consequence} of the words in $\mathcal{R}$ and the group axioms.
\end{definition}

That is, $w$ is the same as the word we get by one or more of the following algebraic transformations : replace any occurrence of $R_i \in \mathcal{R}$ in $w$ by the empty word; and replace any occurrence of $gg^{-1}$ or $g^{-1}g$ in $w$ by the empty word for any $g \in S$.

The elements of $\mathcal{R}$ are called relators. A group presentation may also be written as
$\langle S \mid u_i = v_i, \quad i \in \mathbf{I}  \rangle$ where $u_i, v_i \in S^*$ as before and $\mathbf{I}$ is an indexing set. An equation of the form $u = v$ in $S^*$ is referred to as a \emph{relation}. The relation $u = v$ is equivalent to the relator $u v^{-1}$ as $u = v \iff u v^{-1} = \id$ where $v^{-1}$ is the inverse of $v$ in $S^*$.  It is worth noting at this point that both a relator and a relation can be thought of as an element of $S^* \times S^*$ as $(R_i, \emptyset)$ and $(u, v)$ respectively. We make use of this formulation later in Section~\ref{sec:rws}.

A group $G$ can have many different generating sets and consequently many presentations.

For example, a presentation for the symmetric group $S_n$ with the generating set $S=\{s_i \mid s_i = (i, i+1), 1 \leq i < n\}$ consists of the relations:
\begin{align*}
s_i^2 = \id & \quad \forall 1 \leq i < n \\
s_i s_j = s_j s_i & \quad \mbox{ if } |i-j| > 1 \\
s_i s_{i+1} s_i = s_{i+1} s_i s_{i+1} & \quad 1 \leq i < n-1.
\end{align*}
This is known as the Coxeter presentation~\citep{humphreys1992reflection}. 
In particular for $S_4$, with the generating set $\{s_1, s_2, s_3\}$, we have the following set of relations
\begin{align*}
(R1) & \quad s_1^2 = \id & (R2) & \quad s_2^2 = \id & (R3) & \quad s_3^2 = \id\\
(R4) & \quad s_1 s_2 s_1 = s_2 s_1 s_2  & (R5) & \quad s_2 s_3 s_2 = s_3 s_2 s_3 
 & (R6) & \quad s_1 s_3 = s_3 s_1. 
\end{align*}

The word $w = [s_2, s_3, s_2, s_1, s_3, s_1, s_2, s_3]$ satisfies  
$s_2 s_3 s_2 s_1 s_3 s_1 s_2 s_3 
=
\id$ in $S_4$, meaning $w 
=
\id$ is an algebraic consequence of the group axioms and the relations in the presentation of $S_4$. This can be seen by rewriting $w$ using the relations in the presentation and the group axioms, for example as follows.
\begin{align*}
s_2 s_3 s_2 s_1 s_3 s_1 s_2 s_3 &= s_2 s_3 s_2 s_1 s_3 s_1 (s_3 s_3) s_2 s_3 & [R3, g\id=g] \\
&= s_2 s_3 s_2 (s_1 s_3 s_1 s_3) s_3 s_2 s_3 & [R6, g\id=g]\\
&= s_2 s_3 s_2 s_3 s_2 s_3 & [R5, R2 \text{ and }R3]\\
&= \id.
\end{align*}

The above example suggests how the relations might be developed into a set of rewriting rules and the process of rewriting carried out in a systematic manner.  
In Section~\ref{sec:rws}, we will formalise the notion of such a rewriting system and discuss the properties that make a rewriting system effective.

\section{Words on a Cayley Graph}\label{sec:cg}

Another useful way to understand relations and rewriting of words in groups is through a Cayley graph. For a group $G$ and a generating set $S \subseteq G$, the Cayley graph $\cayley(G,S)$ of $G$ with respect to $S$ is a directed graph that has a vertex for each element of $G$.  There is an edge between vertices $g$ and $h$ if $gh^{-1} \in S$. That is, there is an edge labeled $s$ between $g$ and $h$ if there is some $s \in S$ such that $sh = g$. 

The labels on the edges in a path from vertex $h$ to $g$ in $\cayley(G,S)$ give a word $w$ in $S^*$ such that $\Gamma(w) = gh^{-1}$. Recall that $\Gamma$ maps a word in $S^*$ to an element in $G$. The length of a group element $g$ is the length of a shortest path between the identity vertex $\id$ and $g$. If the edges of this graph are assigned weights, we can talk about the weighted path length between two vertices. 
In particular, a path in $\cayley(G,S)$ from a vertex $g$ to itself concatenates to give a word $w$ that represents $\id$.
Thus, relators from the group are represented by closed paths (loops) in the Cayley graph.
Since the Cayley graph of a group is vertex-transitive, any node in the graph may be fixed as the identity vertex.

In the context of genome rearrangement models, a permutation is a genome arrangement. The generating set is the set of allowed rearrangements under this model. For instance, when inversions are considered to be the only allowed rearrangements, then the generating set is $S=\invset$. The set of all genome arrangements is the genome space which corresponds to the vertex set of the Cayley graph $\cayley(G, \invset)$.

The process of rewriting words using relators is equivalent to deforming a path in a Cayley graph using loops to identify `shortcuts'. As we have seen, a word in $S^*$ that equals $\id$ can be rewritten using the relators in the group presentation and the group axioms. On a Cayley graph, this can be understood as a closed path being constructed using the closed paths in the presentation as building blocks.

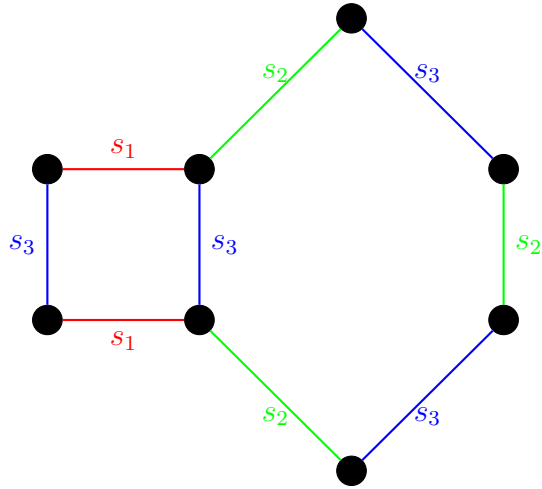
\begin{figure}[ht]
\begin{center}
  \begin{tikzpicture}[scale=2]
  \tikzset{mystyle/.style={draw,circle,fill, black}}  
  \node[mystyle] (e) at (0,0) {};
  \node[mystyle] (s1) at (1,0) {};
  \node[mystyle] (s1s3) at (1,1) {};
  \node[mystyle] (s3) at (0, 1) {};
  \draw[red, thick] (e) -- node[midway,below] {$s_1$} (s1);
  \draw[blue, thick] (e) -- node[midway,left] {$s_3$} (s3);
  \draw[blue, thick] (s1) -- node[midway,right] {$s_3$} (s1s3);
  \draw[red, thick] (s3) -- node[midway,above] {$s_1$} (s1s3);
  \node[mystyle] (s3s2) at (2,2) {};
  \node[mystyle] (s3s2s3) at (3,1) {};
  \node[mystyle] (s3s2s3s2)  at (3,0) {};
  \node[mystyle] (s3s2s3s2s3) at (2, -1) {};
  \draw[thick, green] (s1s3) -- node[midway,above] {$s_2$} (s3s2);
  \draw[thick, blue]  (s3s2) -- node[midway,above] {$s_3$} (s3s2s3);
  \draw[thick, green] (s3s2s3) -- node[midway,right] {$s_2$} (s3s2s3s2);
  \draw[thick, blue] (s3s2s3s2) -- node[midway,below] {$s_3$} (s3s2s3s2s3); 
  \draw[thick, green] (s3s2s3s2s3) -- node[midway,below] {$s_2$} (s1);
  \end{tikzpicture}
  \caption[Loops on a Cayley Graph]{A word in $S^*$ that equals $\id$ can be rewritten using the relators in the group presentation and the group axioms. For instance, the word $s_1s_2s_3s_2s_3s_2s_1s_3 = \id $ in $S_4$.
  Walking along an edge is equivalent to multiplying by a generator (edge label). Starting in the top left corner and tracing the word clockwise, we get the word $s_1s_2s_3s_2s_3s_2s_1s_3$. The closed path  $s_1s_2s_3s_2s_3s_2s_1s_3$ is constructed using the relators $s_1s_3s_1s_3$ and $s_2s_3s_2s_3s_2s_3$.}
  \label{fig:cayleyg}
  \end{center}
  \end{figure}
\section{Rewriting Systems}\label{sec:rws}
As discussed earlier, a set of relators $\mathcal{R}$ is a subset of $S^* \times S^*$ and thus, $\mathcal{R}$ defines a binary relation on $S^*$.  We write $\rws^*$ for the reflexive transitive closure of $\rws$. 
This relation is made compatible with the multiplication in $S^*$ as follows: 
\begin{quote}
if $(l, r) \in \mathcal{R}$, then for words $u,v$ in $S^*$, we say that\\ 
$ulv$ \emph{rewrites to} $urv$.
\end{quote}

If we impose the constraint that $\mathcal{R}$ be antisymmetric, (i.e., $(l,r) \in \mathcal{R} \implies (r,l) \notin \mathcal{R}$), then the reduction process becomes directional. In this case, $\mathcal{R}$ is referred to as a \emph{rewriting system}.  For $(l,r) \in \rws$, we will write $l \writesto r$ and refer to $l$ as the left side of the rule and $r$ as its right side.  If $(x, y) \in \rws^*$, this means that $x$ can be reduced to (rewritten as) $y$ using the rules in $\rws$. We will write this as $x \writestoc{\rws^*} y$.

A word $w \in S^*$ is said to be \emph{reducible} with respect to $\rws$ if there is some $z \in S^*$ such that $w \writestoc{\rws} z$. If no such $z$ exists, then $w$ is said to be \emph{irreducible} with respect to $\rws$.

In applying the rewriting rules to rewrite a word, one may have to make choices at each step. For instance, a word may contain the left sides of more than one rule in $\rws$. For the process of rewriting to be effective, we need to ensure that a given word can be reduced to a unique irreducible word. In addition to this, an essential requirement is that this irreducible representative can be obtained by the application of rewriting rules in a finite number of steps. Formally, we talk about \emph{confluence} and \emph{termination} of a rewriting system.

	\begin{figure}[t]
	\centering
	\begin{tikzcd} [scale=2]
	                        &  u \arrow{dl}[above]{*} \arrow[dr, "*"] & \\
	v \arrow{dr}[below]{*} &                                  & w \arrow[dl, "*"] \\
	& x &
	\end{tikzcd}
	\caption[Confluence]{In a confluent rewriting system, if a word $u$ can be reduced to the words $v$ and $w$, then $v$ and $w$ can be reduced to some word $x$.}
	\label{fig:confluence}
	\end{figure}
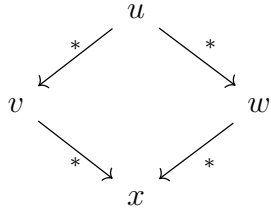

\begin{definition}[Terminating rewriting system]
A rewriting system $\rws$ over $S^*$ is said to be \emph{terminating} if there is no infinite sequence of words $w_i \in S^*$ such that $w_0 \writesto w_1 \writesto \dots  w_k \dots $.
\end{definition}

\begin{definition}[Confluence]\label{d:confluence}
A rewriting system $\rws$ over $S^*$ is said to be \emph{confluent}  if 
 for all  $u,v, w \in S^*$, \text{ if }
$u \writestoc{\rws^*} v$ \text{ and } 
$u \writestoc{\rws^*} w$,
{ then there exists }
$ x \in S^*$ \text{ such that } $v \writestoc{\rws^*} x \mbox{ and } w \writestoc{\rws^*} x.$
(See Figure~\ref{fig:confluence}).
\end{definition}

A set of defining relations (or relators) in a presentation can be turned into a rewriting system. To ensure that the rewriting system thus created is terminating and confluent, we will need to do some more work.

\subsection{Termination}

The termination of a rewriting system $\rws$ can be established by imposing a reduction order on the set $S^*$. A \emph{reduction order} on $S^*$ is a transitive relation $ > $ such that for any $s, t \in S^*$, 
\begin{itemize}
\item exactly one of the following holds: $s > t$, $t > s$ or $s=t$, and
\item $s > t \implies asb > atb, \text{ for all } a, b \in S^*$, and
\item there is no infinite sequence of elements $s_0, s_1, \dots, s_i, \dots$ of $S$ such that
$s_i > s_{i+1}$. 
\end{itemize}

The idea behind imposing $>$ on $S^*$ is that if $u \writesto v \implies u > v$, then an infinite sequence of words $w_i$ such that $w_i \writesto w_{i+1}$ induces an infinite decreasing sequence under $>$. Since the latter is not possible, $\rws$ must terminate.

We now define a reduction order on $S^*$, using a weight function on $S$.

\begin{definition}[Weighted Lexicographic Order]
Let $S$ be a non-empty finite set. Fix any ordering $\succ$ on the elements of $S$. Let $\omega : S \rightarrow \mathbb{R}^+$ be a function that assigns a positive weight to each element of $S$. Let $u=s_1s_2\cdots s_k$ and $v=t_1t_2\cdots t_l$ be in $S^*$. Define $u > v$ if either
\begin{enumerate}
	\item $\sum_{i=1}^k{\omega(s_i)} > \sum_{i=1}^{l}{\omega(t_i)}$, or
	\item $\sum_{i=1}^k{\omega(s_i)} =\sum_{i=1}^{l}{\omega(t_i)} \mbox{ and } s_i \succ t_i \mbox{ where } i=min\{p : s_p \neq t_p\}$.
\end{enumerate}
\end{definition}

It is easy to see that the weighted lexicographic order is a reduction order. 
\begin{prop}
For any finite set $S$, weighted lexicographic order is a reduction order on $S^*$.
\end{prop}

\subsection{Confluence}
If certain mathematical conditions are satisfied, a rewriting system can be transformed into a confluent rewriting system through a procedure due to \citet{knuth1983simple}. We will discuss the Knuth-Bendix algorithm and the properties necessary for it to return a terminating, confluent rewriting system later in this section after introducing the necessary definitions.

\begin{definition}[Critical Pair]
Let $\rws$ be a rewriting system over $S^*$.
Let $u_1 a \writesto v_1$ and $a u_2 \writesto v_2$ be two rules in $\rws$ where $u_i, v_i, a \in S^*, a \neq \id$. That is, a non-empty suffix of the left hand side of a rule overlaps a prefix of the left hand side of another rule. Rules $u_1 a \writesto v_1$  and $a u_2 \writesto v_2$ are said to \emph{overlap}. The word $w = u_1au_2$ can be reduced to both $v_1u_2$ and $u_1 v_2$. The words $v_1u_2$ and $u_1 v_2$ are said to constitute a \emph{critical pair}.
\end{definition}

A critical pair $(v_1u_2, u_1 v_2)$ is said to be \emph{resolved} if there exists $w \in S^*$ such that $v_1u_2 \writestoc{\rws^*} w$ and $u_1 v_2 \writestoc{\rws^*} w$. 

\begin{theorem}[{\citep[Lemma 2.7.2]{baader1999term}}]
\label{thm:cpl}
A terminating rewriting system is confluent if and only if all its critical pairs are resolved.
\end{theorem}

The power of Theorem~\ref{thm:cpl} derives from the fact that it allows us to ascertain the (global) confluence of a rewriting system by checking for confluence locally. This suggests a simple procedure for making a rewriting system confluent. Resolve each critical pair $(u, v)$ by adding a rule $u \writesto v$ if $u > v$ and $v \writesto u$ otherwise. This is the gist of the Knuth-Bendix algorithm. However, we still need to ensure that this loop of adding a rule and checking if any critical pairs remain to be resolved will terminate. In fact, the Knuth-Bendix algorithm is guaranteed to terminate with a confluent and terminating rewriting system if the equivalence relation generated by $\rws$ has finitely many equivalence classes	 \citep[Corollary 12.21]{holt2005handbook}. 

\begin{definition}[Local Confluence](\cite{baader1999term})
A rewriting system $\rws$ over $S^*$ is said to be \emph{locally confluent}  if 
 for all  $u,v, w \in S^*$, \text{ if }
$u \writestoc{\rws} v$ \text{ and } 
$u \writestoc{\rws} w$,
{ then there exists }
$ x \in S^*$ \text{ such that } $v \writestoc{\rws} x \mbox{ and } w \writestoc{\rws} x.$
\end{definition}

Note that local confluence differs from confluence (Definition~\ref{d:confluence}) in that here, relations are from $\rws$ rather than its reflexive transitive closure $\rws^*$.

For a terminating and locally confluent rewriting system,  each equivalence class under the closure of the relation generated by $\writesto$ contains a unique, irreducible element.  
Since each element of $S^*$ maps to a group element, each unique, irreducible element maps to a unique group element. The number of equivalence classes in $S^*$ must be finite if the group generated by $S$ is finite. Thus for a finite group, the Knuth-Bendix algorithm will give us the requisite set of rewriting rules. 

The upshot of this observation is that for a genome rearrangement model, where the rearrangement operators are invertible, the Knuth-Bendix algorithm is guaranteed to generate a finite, confluent, terminating  rewriting system since we are dealing with finite groups. The restriction that the operators be invertible is necessary to ensure that the operators generate a group. 

In Section~\ref{sec:examples}, we construct rewriting systems for two different weighted inversion models and use them to find the weighted distance for genomes. 

\section{Implementation and Biological Examples}
\label{sec:examples}
The first model consists of unsigned permutations on a linear genome with $7$ regions. The set of inversions $\invset$ consists of all inversions $t_{i,j}$, for $1 \leq i < j \leq 7$, as defined in Section~\ref{subs:invsys}. This set generates the symmetric group $S_7$. A simple monotonic weighting function $\omega$ is given by $\omega(t_{i,j}) := j - i$.
For $1\leq i_1 < j_1\leq 7$ and  $1\leq i_2 < j_2\leq 7$, let $m=m(i_1, j_1, i_2, j_2)$ be the smallest non-negative integer such that $(t_{i_1,j_1} t_{i_2,j_2})^m = 1$.
 A presentation for the group with this generating set is
\[
G = \{ \invset \mid (t_{i_1,j_1}t_{i_2,j_2})^{m},\ m=m(i_1, j_1, i_2, j_2),\ 1\le i_1 < j_1 \le 7,\ 1\le i_2 < j_2 \le 7 \}.
\]

We used the software package KBMAG \citep{holtkbmag} to run Knuth-Bendix on this presentation. The resulting confluent rewriting system consists of 6220 rules. KBMAG can also use the rewriting system to find a minimal representative for a given group element.  The weighted distance for the group element is then defined to be the weight of the unique minimal representative.

\begin{figure}[t]
	\centering
	 \begin{tikzpicture}[scale=0.75]
	  \draw [ultra thick] (-1.5, 2) node [left] {A} -- (0,0) -- (0,-1) -- (-3, -2) node [below] {B} ;
	  \draw [ultra thick] (1,2) node [right] {C} -- (0,0) -- (0, -1) -- (3,-3) node
	  [right] {D};
	  \node [left] at (-0.75, 1) {1.5};
	  \node [right] at (0, -0.5) {0.5};
	  \node [below] at (-1.5, -1.5) {2.5};
	  \node [right] at (0.5, 1) {5.5};
	  \node [right] at (1.7, -2) {3.5};
	  \node at (.3,-4) {Weighted Distances};
	 \end{tikzpicture}
	 \hspace{3mm}
	\centering
	 \begin{tikzpicture}[scale=0.75]
	  \draw [ultra thick] (-1.5, 2) node [left] {A} -- (0,0) -- (0,-1) -- (-2.5, -2.5) node [below] {C} ;
	  \draw [ultra thick] (2,2) node [right] {B} -- (0,0);
	  \draw [ultra thick] (0, -1) -- (3,-3) node [right] {D};
	  
	  \node [left] at (-0.75, 1) {2.75};
	  \node [left] at (0, -0.5) {0.75};
	  \node [below] at (-1.1, -1.75) {2.25};
	  \node [right] at (1, 1) {1.25};
	  \node [right] at (1.7, -2) {3.75};
	  \node at (.3,-4) {GRIMM Distance};
	 \end{tikzpicture}
		 \hspace{3mm}
	\centering
	 \begin{tikzpicture}[scale=0.75]
	  \draw [ultra thick] (-1.5, 2) node [left] {A} -- (0,0) -- (0,-1) -- (-2.5, -2.5) node [below] {D} ;
	  \draw [ultra thick] (2,5) node [right] {C} -- (0,0) -- (0, -1) -- (3,-3) node
	  [right] {B};
	  \node [left] at (-0.85, 1) {2};
	  \node [left] at (0, -0.5) {1};
	  \node [below] at (-1.45, -1.25) {4};
	  \node [right] at (1, 2.5) {7};
	  \node [right] at (1.7, -2) {4};
	  \node at (.3,-4) {Coxeter Distances};
	 \end{tikzpicture}
	\caption[A phylogeny drawn using the weighted distances on linear genomes]{The topologies produced with the distance matrices from the different distance algorithms as input. The four permutations are $ A = (1,4)(3,7,6)$, $B = (1,3,7,5,2,6,4)$, $C = (3,4,6)$ and $D = (1,7,6,4,2,3,5)$.}\label{fig:linearphlyo}
	\label{fig:1wted} 
	\label{fig:1grimm} 
	\label{fig:1cox} 

	\end{figure}
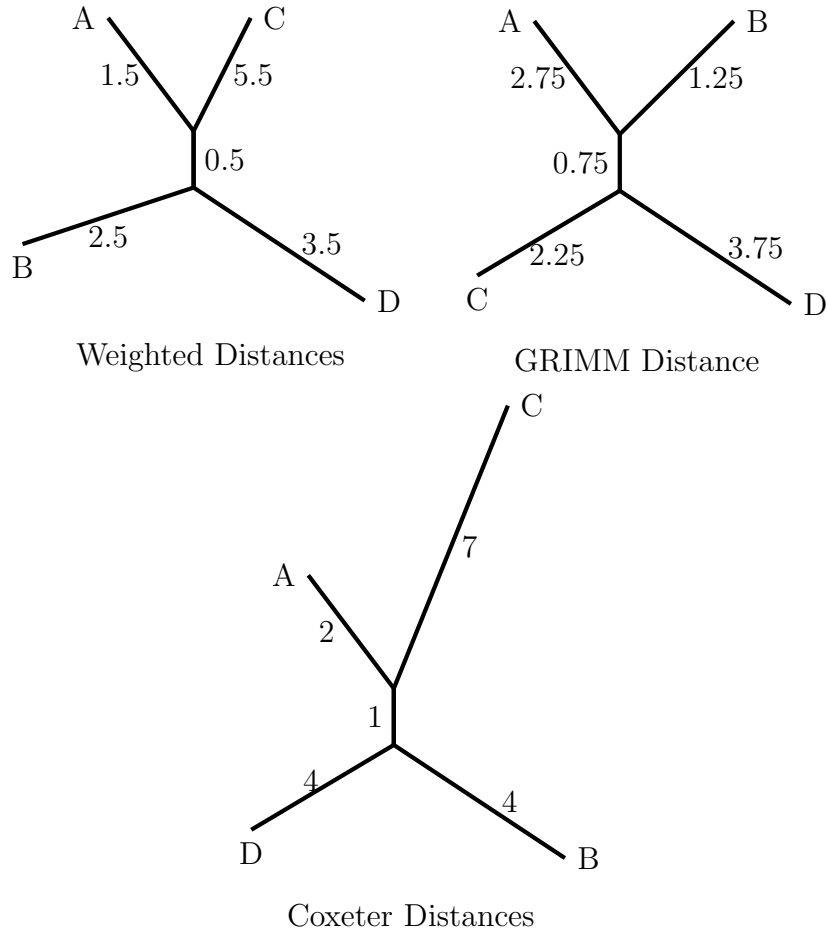

We generated $4$ random permutations in $S_7$ and determined the weighted distance matrix using the weight function $\omega$, which was fed into RPhylip \citep{revell2014rphylip,felsenstein1993phylip} to construct a phylogenetic tree (Figure~\ref{fig:1wted}). For the same permutations, we also constructed phylogenies with the distance matrix from GRIMM \citep{tesler2002grimm} and the Coxeter distance matrix (also Figure~\ref{fig:1grimm}). GRIMM assigns unit weight to all inversions. Coxeter generators are reversals of adjacent regions $(i, i+1)$, and so the inversion model underlying Coxeter distance assigns unit weight to reversals of adjacent regions and infinite weight to all other reversals.

The three topologies presented in Figure~\ref{fig:linearphlyo} differ from each other in either the clustering of nodes or the branch lengths. An important point to be noted is that the weighted distance model results in the clustering $AC|BD$ while the uniform weight model (GRIMM) clusters the nodes as $AB|CD$. This difference is interesting since both the methods have the same set of inversions but different weights asisgned to the generators.

Our second example deals with circular rather than linear genomes.  In this case, we will use the method for returning a minimal distance for adjacent inversions~\citep{egrinagy2013group} to confirm the weighted distance methods in this paper are effective.

To construct the rewriting system, we used the same set of generators and relations as those in the circular inversion model presented in \citet{egrinagy2013group}. The generating set consists of the inversions of adjacent regions $(i, i+1)$ for $1 \leq i < n$ and the inversion $(1,n)$ that allows swap the positions $n$ and $1$. Following the notation in \citet{egrinagy2013group}, we denote these generators by $s_i$. The generating set is thus $\{s_i \mid i=1,2,\dots n\}$ and the relations are:

\begin{align*}
&s_i^2 = \id & \text{ for each } i=1, \dots n, \\
&s_i s_j = s_j s_i & \text{ if } i-j \bmod n \neq 1, \\
&s_i s_{i+1} s_i = s_{i+1} s_i s_{i+1} & \text{ for each } i=1, \dots n-1, and\\
&s_n = s_{n-1} s_{n-2} \dots s_2 s_1 s_2 \cdots s_{n-2} s_{n-1}.
\end{align*}

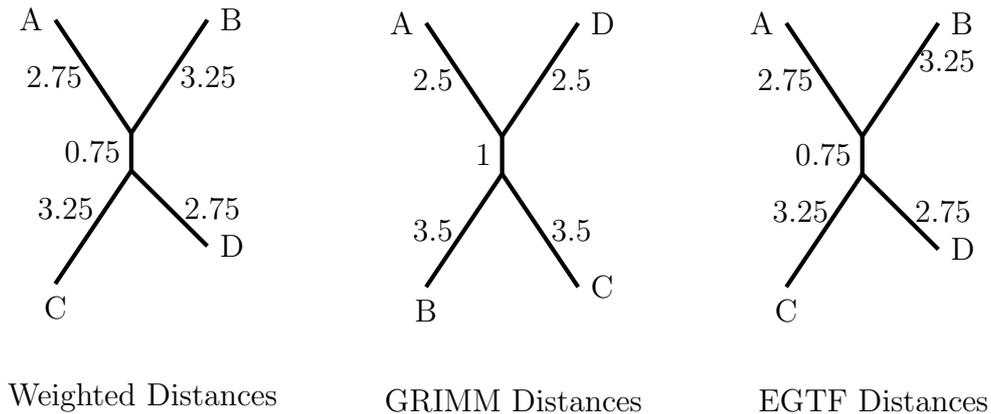
\begin{figure}[t]
	\centering
	 \begin{tikzpicture}[scale=0.5]
	  \draw [ultra thick] (-2, 3) node [left] {A} -- (0,0) -- (0, -1) -- (-2, -4) node [below] {C} ;
	  \draw [ultra thick] (2, 3) node [right] {B} -- (0,0);
	  \draw [ultra thick] (0, -1) -- (2,-3) node [right] {D};
	  \node [left] at (-1, 1.5) {2.75};
	  \node [left] at (0, -0.5) {0.75};
	  \node [left] at (-0.7,-2) {3.25};
	  \node [right] at (1,1.5) {3.25};
	  \node [right] at (1.1, -2) {2.75};
	  \node at (.3,-7) {Weighted Distances};
	 \end{tikzpicture}\hspace{1cm}
	\centering
	 \begin{tikzpicture}[scale=0.5]
	  \draw [ultra thick] (-2, 3) node [left] {A} -- (0,0) -- (0, -1) -- (-2, -4) node [below] {B} ;
	  \draw [ultra thick] (2, 3) node [right] {D} -- (0,0);
	  \draw [ultra thick] (0, -1) -- (2,-4) node [right] {C};
	  \node [left] at (-1, 1.5) {2.5};
	  \node [left] at (0, -0.5) {1};
	  \node [left] at (-1,-2.5) {3.5};
	  \node [right] at (1,1.5) {2.5};
	  \node [right] at (1, -2.5) {3.5};
	  \node at (.3,-7) {GRIMM Distances};
	 \end{tikzpicture}\hspace{1cm}
	\centering
	 \begin{tikzpicture}[scale=0.5]
	  \draw [ultra thick] (-2, 3) node [left] {A} -- (0,0) -- (0, -1) -- (-2, -4) node [below] {C} ;
	  \draw [ultra thick] (2, 3) node [right] {B} -- (0,0);
	  \draw [ultra thick] (0, -1) -- (2,-3) node [right] {D};
	  \node [left] at (-1, 1.5) {2.75};
	  \node [left] at (0, -0.5) {0.75};
	  \node [left] at (-0.6,-2) {3.25};
	  \node [right] at (1.2,2) {3.25};
	  \node [right] at (1.1, -2) {2.75};
	  \node at (.3,-7) {EGTF Distances};
	 \end{tikzpicture}
    \caption[A phylogeny drawn using the weighted distances on circular genomes]{The topologies produced for circular genomes on 8 regions with the distance matrices from the different distance algorithms. The four permutations are $A = (2,3)(6,8)$, $B = (1,7,6,8)(4,5)$, $C = (1,5,6,4,3,8,7)$ and $D=(1,2,4)(5,6,7)$.}\label{fig:circphlyo}
    \label{fig:circwted} 
    \label{fig:circgrimm} 
	\label{fig:circegtf} 
	\end{figure}

All the generators are assigned unit weight as in the circular inversion model of \citet{egrinagy2013group}, which we will refer to as EGTF. 
We use this presentation as an input to KBMAG.  The confluent rewriting system in this case has $6622$ rules.  Once again, we generated $4$ random permutations in the group and found the distances using KBMAG, GRIMM, and EGTF. The latter method has been implemented by the authors in the package BioGAP \citep{egri2014bacterial} for the software GAP \citep{Sch97}. Since EGTF and the method in this paper have the same generating set and the same weights, the phylogenies produced using the distance matrix from the rewriting system and that from the EGTF method should be identical, as indeed they are.
The resulting phylogenies produced using RPhylip are presented in Figure~\ref{fig:circphlyo}. 

The topologies produced by the distances derived from the rewriting system and EGTF model are the same as expected since both these methods give the exact minimal reversal distance between two circular permutations. Both the methods have been set up to factor in the rotational and reflection symmetry of a circular genome.

\section{Discussion and Future Work}
\label{sec:future}
Researchers have recognised the need for methods to determine weighted distances in the field of genome rearrangement right from the start. Beginning with the pioneering work of \citet{sankoff1992edit} and \citet{sankoff1992gene}, a number of approaches have been tried. While they differ in the techniques employed, a common feature of the previous studies is that the proposed algorithms are  tied to a particular model of rearrangement. The novelty of our work is that the framework presented here can be adapted to a wide variety of models. In addition, to the best of our knowledge, this work presents the first use of the theory of rewriting systems to a problem in comparative phylogenetics. 

The current approach however has some limitations, which present opportunities for interesting research.  For instance, the method presented here can only be used  with invertible rearrangement operators. The use of other algebraic structures such as a semigroup might allow this restriction to be removed allowing more rearrangement models to be included.

Another important limitation is that the method works by distorting an existing path (in terms of the operators in the model) between two genomes into an optimal path.  This is not a problem if the model includes all inversions, or all adjacent inversions --- in which case methods such as GRIMM and EGTF can provide an initial path. 
However, for some models, finding a path between two arbitrary genomes may be non-trivial. 

\begin{table}[h!]
  \centering
  \begin{tabular}{lr}
  \toprule
  $n$ &   Rules \\
  \midrule
  3 & 9 \\
  4 & 44  \\
  5 & 204 \\
  6 & 1049 \\
  7 & 6220 \\
  \end{tabular}
  \caption[Number of rewriting rules]{The size of a confluent rewriting system grows very quickly with the number of regions $n$. The rewriting system is for a weighted lexicographic order with weight of $t_{i,j}$ being $j-i$ for a linear genome and generators $(i, i + 1)$ for $i \in \{1, 2, \dots, n - 1\}$. }\label{table:numrules}
\end{table}

Even in the case where such a path is known, for instance in the inversion model, the other deficiency at the moment is the lack of a software implementation. We have used KBMAG to derive the confluent rewriting systems. However, the use of KBMAG with finite groups is not recommended by the authors as it is optimised for infinite groups. It is not surprising therefore that for larger values of $n$, KBMAG cannot return a confluent rewriting system even though we know that a confluent system exists. The size of a confluent rewriting system increases very quickly with $n$ (see Table~\ref{table:numrules}). Thus an efficient software implementation of Knuth-Bendix optimised for finite groups and in particular for models arising from biology would be very useful.

The application of rewriting systems to a new problem also gives rise to new mathematical questions. For instance, it would be interesting to investigate the effect of the weighting function used on the size and efficiency of the rewriting system.

\bibliographystyle{plainnat}

\end{document}